\documentclass[12pt,oneside]{amsart}
\usepackage{verbatim,amsfonts,color}
\usepackage[mathscr]{eucal}
\usepackage{graphicx}

\theoremstyle{plain}
\newtheorem{Thm}{Theorem}
\newtheorem{Cor}{Corollary}
\newtheorem{Lem}{Lemma}
\newtheorem{Prop}{Proposition}

\theoremstyle{definition}
\newtheorem{Deff}{Definition}
\newtheorem{Rmk}{Remark}

\newcommand{\mat}[4]{\left( \begin{array}{cc}
{#1} & {#2} \\
{#3} & {#4} \end{array} \right) }

\newcommand{\threesurf}{\ensuremath{\Gamma^3 \backslash {\mathcal H}} }

\usepackage{hyperref}

\begin{document}


\title[Lagrange value 3 gives Hausdorff  dimension zero]{Hausdorff dimension of the set of real numbers of Lagrange value three}
\author{Thomas A. Schmidt}
\address{Oregon State University\\ Corvallis, OR 97331}
\email{toms@math.orst.edu}
\author{Mark Sheingorn}
\address{Enfield, NH 03748} 
\email{marksh@alum.dartmouth.org}
\date{19 January, 2012}



\begin{abstract}   We  show that the set of real numbers of Lagrange value 3 has Hausdorff dimension zero.   In fact,  we prove the appropriate generalization of this statement for each element of the Teichm\"uller space of the commutator subgroup of the classical modular group.     \\
We also show that for each positive integer $m$, about every rational number there is an open interval of transcendental endpoints that is free of reals of Lagrange value greater than $m$.
\end{abstract}

\maketitle

\section{Introduction}
The quality of approximation by rationals to an irrational number $\xi$ is measured by its Lagrange value 
\[ 
\mu(\xi) = \sup \{ \,h \mid \, \vert \xi - p/q \vert < 1/h q^2 \, \mbox{for infinitely many}\, p/q \in \mathbb Q \} \,.
\]
The Lagrange spectrum is the set of all values $\mu(\xi)$, this spectrum  begins with a discrete sequence with lowest value $\sqrt{5}$ and whose limit point is $3$.    For each of the values below three there are countably many real numbers giving each value; however,  there are uncountably many reals with Lagrange value $3$.    These results were shown by A.~A.~Markoff at the end of the 19th century, in slightly different terms.   (Indeed, the Lagrange spectrum agrees with the ``Markoff spectrum'' up to a value well  --- for our purposes ---  above three.)  See \cite{CF} for these various facts and proofs for some of the following statements.    

A line of argument that goes back at least to L.~E.~Ford early in the 20th century,  lead to the following notion.    Suppose that $\Gamma$ is a zonal  Fuchsian group (recall that this means that $\infty$ is a parabolic fixed point of $\Gamma$),   and $\xi$ is an element of the limit set of $\Gamma$ that is not a parabolic fixed point.   Then the Lagrange value of $\xi$ with respect to $\Gamma$ is 
\[ 
\begin{aligned}
\mu_{\Gamma}(\xi) = \sup \{\,  h \mid \, & \mbox{there exists an infinite sequence}\, \{V_j\} \subset \Gamma\\
\, &\;\;\;\mbox{such that}\, \vert V_j(\infty) - V_j(\xi)\vert\, \mbox{converges to} \, h \} \,.
\end{aligned}
\]
 In particular, the diameter of any $\Gamma$-image of the vertical line ending at $\xi$ is less than $\mu_{\Gamma}(\xi)$.      When $\Gamma$ is $\text{PSL}_2(\mathbb Z)$,  the Lagrange value of $\xi$ with respect to $\Gamma$ is exactly the classical Lagrange value.      The Lagrange value of a real also equals its Lagrange value with respect to $\Gamma^3$, the subgroup  generated by the cubes of $\text{PSL}_2(\mathbb Z)$ --- this, since left cosets of this subgroup have representatives that are translations.       We thus in fact prove that the set of reals of $\Gamma^3$-Lagrange value $3$ has Hausdorff dimension zero.   With appropriate phrasing, our proof works for all groups in the Teichm\"uller space of $\Gamma^3$.\\

Our use of $\Gamma^3$ is far from unmotivated. 
H. Cohn \cite{C} discovered a tight connection between certain  closed geodesics on the once-punctured hyperbolic torus uniformized by the commutator subgroup $\Gamma'$ of the modular group $\text{PSL}_2(\mathbb Z)$ and the Markoff tree of integer triples in terms of which the initial discrete part of the Lagrange spectrum can be given.     A.~Schmidt \cite{Sch} showed how to extend this to any  once-punctured hyperbolic torus.    Lehner and Sheingorn and then Beardon-Lehner-Sheingorn \cite{BLS} related this also to the simple closed geodesics on the surface uniformized by the principal congruence subgroup of level three, $\Gamma(3)$.  A.~Haas \cite{H}  showed in particular that Cohn's relationship involved the {\em simple} closed geodesics of the surface uniformized by $\Gamma'$.  Thereafter,  Sheingorn \cite{Sh}  showed that the relationship holds also for the surface uniformized by $\Gamma^3$.      This is a common supergroup of $\Gamma'$ and $\Gamma(3)$;  the three uniformized ``surfaces'' (in the presence of torsion elements, a group uniformizes rather an orbifold) have their simple closed geodesics in respective one-to-one correspondences.   Using A.~Schmidt's results, Sheingorn also showed that a similar relationship holds throughout the Teichm\"uller spaces of the $\Gamma'$, $\Gamma(3)$ and $\Gamma^3$.   (A.~Haas \cite{H} studied diophantine approximation questions on these various surfaces in terms of penetration of geodesics into a cusp.   The introduction of Haas-Series \cite{HS}  gives a careful discussion of the Lagrange spectrum of points in these terms.)\\

In a previous paper \cite{scgMcId},  we reproved McShane's \cite{Mc} celebrated identity  
\[
\sum_{\gamma}\, \dfrac{1}{1 + e^{\ell(\gamma)}} = \dfrac{1}{2}\;,
\]
 where the sum is taken over all simple closed geodesics of any fixed hyperbolic once-punctured torus, and $\ell(\gamma)$ is the length of the geodesic.

 Our approach was to show that the Hausdorff dimension of a set related to the ``height achieving'' simple geodesics is zero.  In fact, we worked on the sphere formed by taking the quotient of the torus by its elliptic involution;  the corresponding uniformizing group has many elliptic elements of order two.  Each element of order two acts so as to define a ``lifting'' region in which no highest points of geodesics can exist.   We showed that a set, indexed over the simple closed geodesics, of these regions meets the ``fundamental horocycle'' in a union of disjoint intervals, and the complement of these intervals has Hausdorff measure zero.  Elementary bookkeeping then gives the identity.       
 
 The argument both reproves the identity and shows that the set of simple geodesics that have a highest point (which always lie on or below the fundamental horocycle)  has Hausdorff dimension zero.   But,  there also exist simple geodesics that never achieve their highest point.   One need say a bit more so as to find the size of the set of these ``non-achieving'' simple geodesics.     We do this in the present paper.

We replace our   ``lifting regions'' intersecting the fundamental horocycle by the ``shadows'' projected onto the real line by (some of) the horocycles  in the group orbit of the fundamental horocycle.     No element of the shadow can be the foot of a simple non-achieving geodesic.    But,  each simple closed geodesic gives rise to a maximal interval of overlapping shadows,  and the length of each such interval equals the length of the intervals defined by the lifting regions.    Thus,  the argument of the previous paper shows that the complement has Hausdorff dimension zero.  
\bigskip

\begin{Thm}\label{t:main} The set of real numbers of Lagrange value equal to three has Hausdorff dimension  zero.   The analogous statement holds throughout the Teichm\"uller space of  hyperbolic genus zero orbifolds with one cusp and three elliptic fixed points of order two.  
\end{Thm}

\bigskip
We remark that in the classical setting,  every real number of Lagrange value at most three has regular continued fraction expansion with partial quotients in the set $\{1, 2\}$.      The set of all real numbers whose partial quotients take on only these two values has Hausdorff dimension larger than $1/2$, as Good \cite{G} showed in 1941; see Jenkins and Pollicott's \cite{JP}  for more recent results.    The modern approach relies on the thermodynamic formalism,  simplified by the fact that two operators stabilize the set.   This approach has roots including the work of Bumby \cite{Bumby},  who showed that the Hausdorff dimension of the set of reals of Lagrange value at most $\sqrt{689}/8$ is less than $1/2$ (see \cite{CF} for a sketch of the proof).   The Lagrange value three set is not determined by merely two operators,  thus this approach would be awkward in the determination of its Hausdorff dimension.

\section{Background:   Teichm\"uller space of spheres double covered by once-punctured tori}\label{s:background}  

For the ease of the reader,  in this section we reproduce portions of \cite{scgMcId}, with some slight editing.  
Note also that since all matrices considered here are of determinant $1$,   we occasionally denote a fourth entry with a $*$.

\subsection{Fricke's Equation and Explicit Groups}

We are interested in explicit lifts of simple closed geodesics.  For this, we use a variation of A. Schmidt's application \cite{Sch} of work of Fricke.  Suppose that  positive real $a, b, c$ satisfy  the Fricke equation

\begin{equation} a^2 + b^2 + c^2 = a b c\,, \label{FrickeEq}
\end{equation}
the elements 
\[T_{0} :=  \mat{0}{-a/c}{c/a}{0}, \; 
T_{1} := \mat{a/c}{*}{b/a}{-a/c}\,,\;
T_{2} := \mat{a-b/c}{* }{1}{-a +b/c}
\]
(of determinant one) generate a group of signature $(0; 2, 2, 2; \infty)$.   
 Note that 
 \[ T_2\cdot T_1\cdot T_0 = S^{a}: z \mapsto z + a\]
  is the fundamental translation of this group. 
A full set of  orbit representatives under the action of the Teichm\"uller group is given when one takes $2 < a \le b \le c < ab/2$; this can be deduced from \cite{Sch}, see also  \cite{W}; we will always assume that our Fricke triples $(a,b,c)$ satisfy this restriction. Note that the
modular case of \threesurf corresponds to
$a=b=c=3$ and in this case $T_{j}$ is the conjugate of $T_{0}$ by the translation $z \mapsto z + j$.    (Note that our $T_{0}$ is not that of \cite{Sch}.)

\subsection{Fixed Point Triples and   Fundamental Domains}

For ease of presentation, in \cite{scgArt} we restricted to the modular case.  However, as we noted, our arguments extend to the full Teichm\"uller case. 

\begin{Prop}\cite{scgArt}\label{pssiFunDomThm}   Fix a Fricke triple $(a,b,c)$ and the corresponding   signature $(0; 2,2,2; \infty)$-orbifold $\Gamma \backslash \mathbb H$.   Each simple closed geodesic of the orbifold has a highest lift which is the axis of 
$S^{a} E$, where $E \in \Gamma$ is elliptic of order two.   There is a factorization 
of $S^{a} E = GF$ as the product of elliptic elements such that 
a highest lifting segment  of this simple geodesic joins the fixed point $f$ of $F$ 
to the fixed point $g$ of  $G$.     Let $e$ be the fixed point of $E$.   A fundamental domain for $\Gamma$  is  given by the
hexagon of vertices: $\infty$, $e$, $f$, $F(e)$, $g$, $a+e$.   In particular, $\{E, F, G\}$ generates $\Gamma$.
\end{Prop} 

  Given a fixed Fricke triple $(a,b,c)$,  we have the corresponding {\em adjusted Fricke equation}
\begin{equation} x^2 + y^2 + z^2 = a \, x y z\,. \label{AdjFrickeEq}
\end{equation}
Note that when $a=3$ the adjusted Fricke equation is exactly Markoff's equation.  

 Recall that the imaginary part,  $Y$,   of a point $X + i Y$  is its {\em height}.  The factorization $S^a E = G F$ can be used to show that the hyperbolic $GF$, whose axis projects to a simple closed geodesic on the surface, has trace $a z$, where $1/z$ is the height of the fixed point of $E$.   
 One can show that there is such a factorization for every simple closed geodesic, and since the adjusted Fricke equation is satisfied by the traces of appropriate triples of simple hyperbolic elements, one finds the following result. 
   
\begin{Cor}\cite{scgArt}\label{gotTriples}  Let $E, F, G$ be as above.  Then the fixed points of $E, F, G$ have respective  heights  $1/z, 1/y, 1/x$,  whose inverses  give a triple satisfing the adjusted Fricke Equation,  and with $z = \max\{x, y, z\}$.   Furthermore, the simple closed geodesic that lifts to  the axis of $S^{a}E$ has height $r_{a}(z) = \sqrt{a^{2}/4 - 1/z^{2}}$.     
\end{Cor}

\begin{proof}   This follows from Theorems 2 and 3 (and their proofs) of \cite{ssEG}, where we give detailed proofs in the modular case.    The only aspect of the proof given there that does not hold in general is that by use of the  map $w \mapsto -\bar{w}$, in the modular case one can further assume that $y\ge x$.
\end{proof}

\begin{Cor}\cite{scgArt}\label{slideMats}  Let $E, F, G$ be as above.  Then there is a real translation conjugating the triple $E, F, G$ to 
\[
E_{0}= \mat{0}{*}{z}{0}, \;  E_{1}= \mat{x/z}{*}{y}{-x/z},
\; E_{2} = \mat
{ax-y/z}{*}{x}{-ax+y/z}\;.
\]
\end{Cor}

\begin{proof}  This follows  as in the proof of Theorem 2 of \cite{ssEG}. 
\end{proof} 

\bigskip

\noindent
{\bf Convention} For the remainder of the paper, unless otherwise stated, we fix a  Fricke triple $(a,b,c)$,  with $2 < a \le b \le c < ab/2$.   This determines a group  $\Gamma := \langle T_0, T_1, T_2\rangle$ as above.    Note that the fundamental translation length of $\Gamma$ equals $a$.
\bigskip

\subsection{Trees of triples}

\begin{Deff} For any $(E, F, G)$ as in Corollary ~\ref{gotTriples}, we define the following maps to triples of elliptic elements of order two.
\[
\begin{aligned}
\nu: (E,F,G) &\mapsto (FEF, G, S^{a}FS^{-a})\\
\rho: (E,F,G) &\mapsto (FGF, F, S^{a}ES^{-a})\\
\lambda: (E,F,G)  &\mapsto (EFE, E, G)\,.
\end{aligned}
\] 
\end{Deff}

\begin{Prop}\label{theMoves}\cite{scgMcId}  Let  $(E, F, G)$ be as above.  Then each of $\nu(E,F,G), \lambda(E,F,G)$, and $\rho(E,F,G)$ is a generating triple of $\Gamma$.  For each of these triples, the corresponding triple of fixed points gives rise to a solution of the  adjusted Fricke equation, by taking  inverses of  heights.Furthermore, if $z \ge \max\{x,y\}$ then the analogous inequality holds upon applying either of $\lambda$ or $\rho$.
\end{Prop}

\begin{Deff}       For $(E, F, G)$ as above,   let $\mathcal T_{\lambda, \rho}(E, F, G)$ denote the tree formed by applying to the triple all finite compositions (including the identity) of $\lambda$  and $\rho$ to $(E, F, G)$.    Let $\mathcal T^{\nu}_{\lambda, \rho}$ denote the tree formed by joining $\mathcal T_{ \lambda, \rho}(T_{0}, T_{1}, T_{2})$  to \newline 
$\mathcal T_{\lambda, \rho}(\, \nu(T_{0}, T_{1}, T_{2})\, )$  with an edge (labeled by $\nu$).   We call $\mathcal T^{\nu}_{\lambda, \rho}$ the {\em adjusted Fricke tree}.
\end{Deff}
 
\begin{Thm}\label{allScg} \cite{scgMcId} The simple closed geodesics of $\Gamma \backslash \mathbb H$ are exactly the projections of the axes of $T_1 S^{-a}$,  $T_2 S^{-a}$ and the projection of the axis of  $ES^{-a} = FG$ for $(E, F, G)$ a node of the tree 
$\mathcal T^{\nu}_{\lambda, \rho}$. 
 \end{Thm} 

 \bigskip  
    
\section{Shadows of horocycles}

    Our goal is to excise from an interval of length $a$ subintervals that correspond to $\xi$ of Lagrange value greater than $a$.  We will eventually show that the result is a Cantor set of Hausdorff dimension zero.  

\begin{figure}
\begin{center}
\scalebox{0.75}{\includegraphics{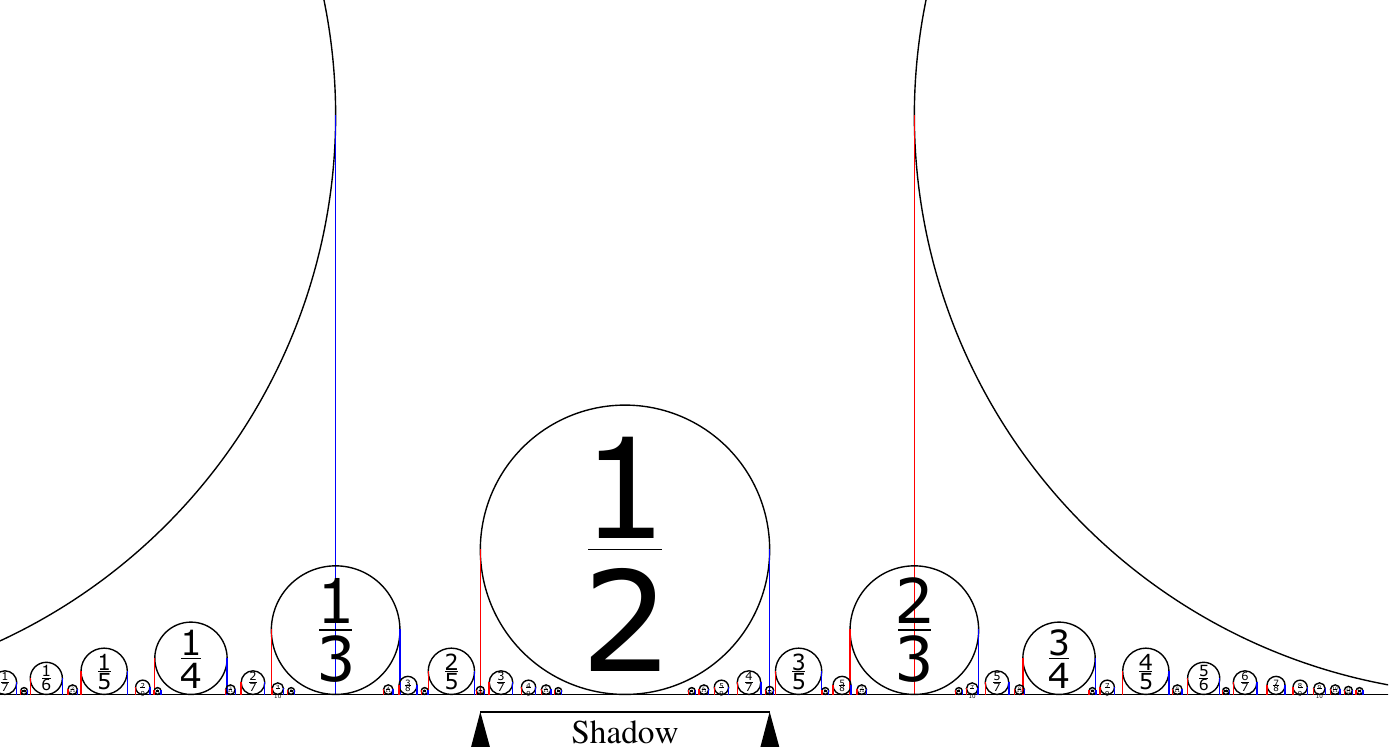}}
\caption{Excision in the unit interval, classical case.  Shadow of horocycle anchored at $1/2$.}
\label{figExcision}
\end{center}
\end{figure}

\bigskip  
\begin{Deff}\label{d:Shadow}    We say that the {\em shadow}  of the   $a$-horocycle anchored at the parabolic point $\alpha/\gamma$ of $\Gamma$   is the closed interval of radius $1/(a \gamma^2)$ centered at $\alpha/\gamma$.        If  $\alpha/\gamma < \beta/\delta$,  we say that $\alpha/\gamma$ and $\beta/\delta$ have {\em overlapping shadows}  if 
\begin{equation}\label{e:overCon}
\alpha/\gamma + 1/(a \gamma^2) >  \beta/\delta - 1/(a \delta^2)\,.
\end{equation}
\end{Deff} 
One can easily discern some some overlapping shadows, in the classical case, in Figure~\ref{figExcision}. 
 
\bigskip  
\begin{Lem}\label{l:theyAreFree}    Any non-parabolic point $\xi$ in the shadow of the   $a$-horocycle anchored at a parabolic point has Lagrange value $\mu_{\Gamma}(\xi)$ at least $a$.   
\end{Lem}
\begin{proof}  This follows from the very definition of  $\mu_{\Gamma}(\xi)$.   
\end{proof}

We prove a key result.

\begin{Prop}\label{p:excisionInter}     Let $M = A S^{-a}$ with $A = \begin{pmatrix}\alpha&\beta\\\gamma&-\alpha \end{pmatrix}\in \Gamma$ an elliptic element of order two such that $M$ is hyperbolic.   Then the interval  
\[  \bigg(\dfrac{\alpha}{\gamma} - \dfrac{a}{2} + \sqrt{\frac{a^2}{4}- 1/\gamma^2}\,, \; \dfrac{\alpha}{\gamma}+ \dfrac{a}{2} -  \sqrt{\frac{a^2}{4}- 1/\gamma^2}\,\bigg)\,\]
is the infinite union of shadows of $M^k \cdot \infty$ with those of $N^k\cdot \infty$ where $N = A S^{a}$.   
\end{Prop} 

The following two lemmas provide the main part of the proof of the proposition.

\begin{Lem}\label{l:markoffOverlapToLeft}    Suppose that  $N = A S^{a}$ is hyperbolic with $A \in \Gamma$ any elliptic element of order two.   Then for all $k \ge 1$,  the parabolic points $N^{k+1}\cdot \infty$ and $N^{k}\cdot \infty$ have overlapping shadows.
\end{Lem}
\begin{proof}    We  have 
\[N = A S^{a} =    \begin{pmatrix} \alpha&\beta - a \alpha\\
                                   \gamma&-\alpha + a \gamma
         \end{pmatrix}\,.
\]                                          
Let $N^k =  \begin{pmatrix} \alpha_k&\beta_k\\
                                   \gamma_k&\delta_k
         \end{pmatrix}$.   
Elementary manipulations easily give that here the overlapping shadows condition \eqref{e:overCon} is equivalent to 
\[ a \gamma_k \gamma_{k+1} (\alpha_{k+1} \gamma_{k} - \gamma_{k+1} \alpha_{k}) + \gamma_{k+1}^{2} + \gamma_{k}^{2}>0\,.\]

Since   $N N^{k} = N^k N$,  we find that 
\[ \alpha_{k+1} = \alpha \alpha_k + (\beta + a \alpha) \gamma_k = \alpha \alpha_k + \gamma \beta_k\]
and
\[ \gamma_{k+1} = \gamma \alpha_k + (-\alpha + a \gamma) \gamma_k = \alpha \gamma_k + \gamma \delta_k\,.\]
This gives in particular that
\[
\begin{aligned}
\alpha_{k+1} \gamma_{k} - \gamma_{k+1} \alpha_{k} &= ( \alpha_k \alpha + \beta_k \gamma)   \gamma_k - (\gamma_k \alpha+  \delta_k \gamma)\alpha_k\\
&= -\gamma\,.
\end{aligned}
\] 
Thus,  our aim is to show that $\gamma_{k+1}^{2} + \gamma_{k}^{2} - a \gamma \gamma_k \gamma_{k+1}$ is positive;  we show in fact that it equals the obviously positive $\gamma^2$.

The identities using $N N^k = N^k N$ also show
\[ \beta_k =  \dfrac{\beta + a \alpha}{\gamma} \gamma_k\,,\]
and 
\[ \alpha_k - \delta_k =  \dfrac{2\alpha - a \gamma}{\gamma} \gamma_k\,.\]
Therefore, 
\[
\begin{aligned}
\gamma_{k+1}^{2} + \gamma_{k}^{2} - a \gamma \gamma_k \gamma_{k+1}&= ( \alpha \gamma_k +  \gamma \delta_k)(\gamma  \alpha_k -   \alpha \gamma_k) + \gamma_{k}^{2}\\
&= \gamma^2 \, [\, \alpha_k \delta_k + \dfrac{\alpha}{\gamma} (\alpha_k - \delta_k + \dfrac{1 -\alpha^2}{\alpha \gamma} \gamma_k) \, \gamma_k\,]\,.
\end{aligned}
\]
That this last equals $\gamma^2$ and is thus clearly positive,  follows from the following. 
\[ \begin{aligned}
\beta_k &=  \dfrac{\beta + a \alpha}{\gamma} \, \gamma_k\\
             &=  \dfrac{1}{\gamma}\bigg(\dfrac{-1-\alpha^2}{\gamma} +  a \alpha \bigg) \, \gamma_k\\
                          &= \dfrac{-1}{\gamma}\bigg(\dfrac{1- \alpha^2 + 2 \alpha^2}{\gamma} - a \alpha \bigg) \, \gamma_k\\
             &=   \dfrac{-\alpha}{\gamma}\bigg(\dfrac{2 \alpha - a   \gamma}{\gamma} + \dfrac{1- \alpha^2 }{\alpha \gamma} \bigg) \, \gamma_k \\
              &=   \dfrac{-\alpha}{\gamma}\bigg(\alpha_k - \delta_k + \dfrac{1- \alpha^2 }{\alpha \gamma} \gamma_k   \bigg)\,.
  \end{aligned}           
\] 
\end{proof}

Similarly, one can establish the following. 

\begin{Lem}\label{l:markoffOverlapToRight}     Suppose that  $M = A S^{-a}$ is hyperbolic with $A \in \Gamma$ any elliptic element of order two.  Then for all $k \ge 1$,  the parabolic points $M^k\cdot \infty$ and $M^{k+1}\cdot \infty$ have overlapping shadows.
\end{Lem}
 
In light of the Proposition ~\ref{p:excisionInter}, we define the the {\em excision interval} of $A = \begin{pmatrix}\alpha&\beta\\\gamma &-\alpha\end{pmatrix}$ as 
\[ I_A =  \bigg(\dfrac{\alpha}{\gamma} - \dfrac{a}{2} + \sqrt{\frac{a^2}{4}- 1/\gamma^2}\,, \; \dfrac{\alpha}{\gamma}+ \dfrac{a}{2} -  \sqrt{\frac{a^2}{4}- 1/\gamma^2}\,\bigg)\,.\]
The {\em width} of excision by $A$ is $w_a(\gamma) = a - \sqrt{a^2 - 4/\gamma^2}$. 
  
\begin{Rmk}   Since no real $\xi$ with $\mu_{\Gamma}(\xi)\le a$ can lie in any of the shadows,   we indeed find that the set of reals with $\mu_{\Gamma}(\xi)\le a$ is contained in the complement of  the excision interval $I_A$.      The proof of our main result, Theorem~\ref{t:main}, relies on the marvelous coincidence that the width of the excision interval $I_A$  is exactly equal to the width of the interval we associated  in \cite{scgMcId} to such an $A$.   In that setting,   $A$ and a certain conjugate of $A$ defined a region where no highest point of a closed geodesic could lie;   the region meets the fundamental horocycle $y=a/2$ in an (excision) interval with this width    $w_a(\gamma)$.
\end{Rmk}

\section{End of proof} 
The set of reals $\xi$ with $\mu_{\Gamma}(\xi)\le a$ is contained in the complement of the union of all excision intervals $I_A$ such that $A \in \Gamma$ is elliptic of order two.
  In fact, arguments in  \cite{scgMcId} allow us to show that the complement to the union of all translates by $x \mapsto x+a$ of the $I_E$ such that  $(E, F, G)$ is a node of the tree $\mathcal T^{\nu}_{\lambda, \rho}$ (and of $I_{T_1}$ and $I_{T_2}$) already results in a set of Hausdorff dimension zero.   The reason for this is that we can form this complement by first excising all translates of $I_{T_1}$ and $I_{T_2}$; following each branch of $\mathcal T^{\nu}_{\lambda, \rho}$ emanating from   $(T_0, T_1, T_2)$,  at the node  $(E,F,G)$  we excise $I_E$ and all of its translates.   One easily verifies that $I_E$ lies between $I_{S^{-a}GS^a}$ and $I_F$.  An easy induction argument using the definition of $\nu, \lambda$ and $\rho$ shows that each of $I_{S^{-a}GS^a}$ and $I_F$ are excised before $I_E$ is (see also Figure~5 of \cite{scgMcId}).      
Applying Corollary~\ref{gotTriples}, we conclude that an interval of length $w_a(z)$ is excised from its ambient interval of length $F\cdot \infty + r_a(y) - (G\cdot\infty - r_a(x)\,)$, and similarly for countably many other intervals.    But, in the proof of Theorem~2 of  \cite{scgMcId} we showed that excising intervals of length  $w_a(z)$ from ambient intervals of  $F\cdot \infty + r_a(y) - (G\cdot\infty - r_a(x)\,)$ results in a Cantor set of Hausdorff dimension zero.

\bigskip 
\section{Further remarks: Shadows in the Classical case}

For any positive integer $m$, the vertical tangents to the $m$-horocycle anchored at a rational $p/q$ have equation $x = p/q \pm 1/(m q^2)$.   Thus  the $m$-horocycles anchored at these new $x$-values cast shadows that certainly overlap with those anchored at $p/q$. 
 Repeatedly applying this observation leads to the following notion.

\begin{Deff} Suppose that $p/q$ is a reduced rational number, then its $m$-{\em successor shadows} are the shadows of the $m$-horocycles anchored at $\sigma_m(p/q) = p/q +1/(m q^2)$ and $\tau_m(p/q)= p/q - 1/(m q^2)$.   The {\em right descendants} of $p/q$ is the collection of the various rationals $\sigma_{m}^{k}(p/q) =  p/q  +1/m q^2 + 1/m(mq^2)^2+ \cdots +m/(m q)^{2^k}= p/q  + m \sum_{n=1}^{k} ( m q)^{-2^n}$.  Similarly, the {\em left descendents} are  the values of the form   $\tau_{m}^{k}(p/q) = p/q  - m \sum_{n=1}^{k} ( m q)^{-2^n}$.
 The $m$-{\em dribble} of $p/q$,  denoted $\mathcal D_{p/q}= \mathcal D_{p/q}(m)$, is the union of the shadows cast by the $m$-horocycles anchored at $p/q$ and at all of its descendants. 
\end{Deff}

\begin{Rmk}
Note that the proof  of Theorem ~\ref{t:main}  does {\em not} use dribbles.   Whereas each successor shadow involved in a dribble is {\em centered} at a vertical tangent line to the horocycle anchored at its predecessor and thus has half of its shadow overlapping with that of its predecessor,   the corresponding ratio goes to zero in the settings of Lemmas~\ref{l:markoffOverlapToLeft} and \ref{l:markoffOverlapToRight}.
\end{Rmk}

 We now describe the dribbles. 

\begin{Prop}\label{p:transcDribble}    For any positive integer $m$ the dribble, $\mathcal D_{p/q}(m)\,$, of a rational number is an open interval with transcendental endpoints.    This interval contains no real number of Lagrange value less than $m$.
\end{Prop} 

\begin{proof}  Since each shadow is a non-empty interval with endpoints being a descendant,  the dribble is the union of overlapping intervals.   The dribble  is thus certainly an interval.  Now,  the $k$th successor to the right of $p/q$ has shadow of right endpoint  $p/q + m\; \sum_{n=1}^{k}\,  (m q)^{-2^n}$,  thus the union has right endpoint $\alpha = p/q + m\; \sum_{n=1}^{\infty} \, (m q)^{- 2^n}$, which is clearly larger than the right endpoint of any of the successor's shadows.    That the left endpoint of the dribble is not included  follows {\em mutatis mutandi}.

  Each of the shadows contains no real of small Lagrange value, thus this holds for their union.  

The endpoints of $\mathcal D_{p/q}$ are given by $\alpha=  p/q + m \; \sum_{k=1}^{\infty} \, (m q)^{-2^k}$ and $\beta=  p/q - m \; \sum_{k=1}^{\infty} \, (m q)^{-2^k}$.   Due to the rationality of $p/q$,  the transcendence of $\alpha$ is equivalent to that of $\beta$.  Thus,  their transcendence is 
implied by  the following result.
\end{proof}

\begin{Prop}\label{p:transcAnyT}    Fix $m\in \mathbb N$,  and the reduced rational number $p/q$.   The real number 
\[\alpha:=  p/q + m \; \sum_{k=1}^{\infty} \, (m q)^{-2^k} \] 
is transcendental.
\end{Prop} 

\begin{proof}  We choose the  regular continued fraction expansion $p/q = [a_0; a_1, \dots, a_n]$ such that $a_n>1$.\\

 \noindent
{\bf Case 1.}  Suppose also that    
\begin{itemize} 
\item  $a_1>1$;
\item $n>1$ is even.
\end{itemize} 

Since $m$ is fixed,  write  $\sigma(p/q) = \sigma_m(p/q)$ and  $\tau(p/q) = \tau_m(p/q)$.         
We show the transcendence of  $\gamma = -a_1 +1/(-a_0+ \lim_{\ell\to \infty}\, \sigma^{\ell}(p/q))$.  

    With van der Poorten \cite{vdP},  let $w' = a_1 \dots a_{n-1}$ and also for any word $v$,   let $\overleftarrow{v}$ denote the word in reversed order. 
We also let $w'' = a_2, \dots, a_{n-1}$;  for any integer $\ell>1$ we let $U(\ell) = U_m(\ell)$ denote the word $\ell, m-1, 1, \ell-1$;   $W = w'' U(a_n) \overleftarrow{w''}$; and $U = U(a_1)$.    We note that the length of $W$ is $2(n-2)+4 = 2n$ and the length of $U$ is 4.

By the Folding Lemma ~\ref{l:folding} below,  we have 
\[ \sigma(p/q) = [a_0;     a_1,  W, a_1]\,,\]
\[  \sigma^2(p/q)  = [a_0; a_1, W, U, \overleftarrow{W}, a_1]\,,\]
\[ \sigma^3(p/q) = [a_0; a_1, W, U, \overleftarrow{W}, U, W, \overleftarrow{U}, \overleftarrow{W}, a_1]\,.\] 
For $k\ge 1$,  let $a_0 a_1 V_{k}' \,a_1$ be the word given by the partial quotients of $\sigma^{k+1}(p/q)$.   By induction,  the partial quotients of $\sigma^{k+2}(p/q) = \sigma([a_0; a_1 V_{k}' \,a_1]) = [a_0; a_1 V_{k}' \,U\, \overleftarrow{V_{k}'}\,a_1]$.

For $k\ge 1$,  let now $V_k$ be the word given by appending $U$ to $V_{k}^{'}$.   The first few of these are thus:    
\[V_1 = W U \overleftarrow{W} U;\;    V_2 = V_1 W \overleftarrow{U} \overleftarrow{W} U;  \; V_3 =V_2 V_1 \overleftarrow{U} \overleftarrow{W} U\,.\]    Induction shows that for $k\ge 3$, 
\[ V_k = V_{k-1} \,\overleftarrow{V_{k-1}'} \, U = V_{k-1} \, V_{k-2}  \, \overleftarrow{V_{k-2}'} \, U\,.\]
The length of $V_k$ is thus clearly twice that of $V_{k-1}$.

Now, for any given $t$, we have that $\gamma = -a_1 +1/(-a_0+ \lim_{\ell\to \infty}\, \sigma^{\ell}(p/q))$ has its initial $t$ partial quotients agreeing with the initial $t$ letters of $V_k$ for any sufficiently large $k$.   The continued fraction expansion of $\gamma$ thus initially agrees with   $V_k V_{k-1}$ for each $k\ge 2$.   That is, the continued fraction expansion agrees with the non-periodic sequence of words $V_k$ of increasing length such that the agreement is with $3/2$ the length of $V_k$ in the sense of Adamczewski and Bugeaud \cite{AB}.  

 To see that $\gamma$ is not quadratic,  it suffices to show that the limit of the $V_k$ is not eventually periodic.   Now,  were this periodic then there would be a $k_0\ge 3$ such that $V_{k_0}$ includes a full period;   since $V_{k_0+1}$ begins with $V_{k_0} V_{k_0-1}$,  already some initial portion of $V_{k_0-1}$ is contained in a period:   Therefore,  the limit of the $V_k$ is purely periodic
and  $V_{k_0}$ is some power of the (minimal) period.   Thus $V_{k_0+1}$ is twice this power, and it follows that $V_{k_0} =  \overleftarrow{V_{k_0}'}U$.   
But, by definition,  $V_{k_0} =  \overleftarrow{V_{k_0}'}U$; therefore,  $V_{k_0}'=\overleftarrow{V_{k_0}'}$.
But, for any $k\ge 2$, 
$V_k = V'_{k-1} \,U\, \overleftarrow{V_{k-1}'} \, U$  and  thus, 
$V'_{k} = V'_{k-1} \,U\, \overleftarrow{V_{k-1}'}$.   Since $U$ is not a palindrome,   there is no $k$ such that $V_{k}'=\overleftarrow{V_{k}'}$.

 By   Theorem 1 of \cite{AB}, $\gamma$ is transcendental.\\   

 \noindent
 {\bf Case 2.}  If   $a_1>1$   but $n>1$ is odd, then the reduced rational
$\tau(p/q) =  [a_0;     a_1,  W, a_1]$, with $W$ as above, is of the type treated in Case 1.   Transcendence of the limit here in turn implies transcendence of $\alpha$ and $\beta$.\\   
 
\noindent 
{\bf Case 3.}  If    $a_1=1$, and  $n>1$ is odd, we find 
\[\tau(p/q) = [a_0;  1, a_2, w''', U(a_n),\overleftarrow{w'''},  a_2+1]\,\]  where   $w''' = a_3\dots a_{n-1}$.     Induction easily confirms that the continued fraction expansion of $\tau^k(p/q)$ is of odd length.   Thus, substituting $W' = w''' U(a_n) \overleftarrow{w'''}$ for $W$ and setting $U = U(a_2+1)$ in the arguments above  again shows the transcendence of a number rationally related to $\alpha$ or $\beta$. \\   

\noindent 
{\bf Case 4.}    If    $a_1=1$, and  $n>1$ is even, then $\sigma(p/q)$ is of the form of $\tau(p/q)$ of the previous case.    Again we find the  transcendence of a number rationally related to $\alpha$ or $\beta$. \\

 \noindent
{\bf Case 5.}  If  $n=1$, then  $a_1>1$ and hence then $\tau(p/q) = [a_0;  U(a_1)]$, and the arguments of Case 1 can be applied. 
     
\end{proof}

\bigskip 
A.~van der Poorten, see the second section of  \cite{vdP},  gave the following version of the Folding Lemma. For the history of this lemma,  see section 10 of \cite{AA} and the notes for section 6.5 of \cite{AS}. 

\begin{Lem}[Folding Lemma]\label{l:folding}    Suppose that the continued fraction expansion is $p/q$ is 
$[a_0; a_1, \dots, a_n]$ in the standard notation, including insisting that $a_n>1$.   If  $m \ge 2$ is a natural number, then 
\[
\begin{aligned} 
 p/q  + (-1)^n/(m q^2) &= [a_0; a_1, \dots, a_n, -m, -a_n, -a_{n-1}, \dots, -a_1] \\
                            &=  [a_0; a_1, \dots, a_n , m-1, 1, a_n-1, a_{n-1}, \dots, a_1]\\
\end{aligned}
\]                            
\end{Lem}


\begin{thebibliography}{Thurst88}

\bibitem{AA} B.~Adamczewski and J.-P.~Allouche,
{\em Reversals and palindromes in continued fractions},   Theoretical Computer Science 380 (2007) 220--237.


\bibitem{AB} B.~Adamczewski and Y.~Bugeaud,
{\em On the complexity of algebraic numbers, {II}. {C}ontinued
              fractions},   Acta Math., 195 (2005), 1--20.
              
\bibitem{AS} J.-P.~Allouche and J. Shallit,
{\em Automatic sequences:
theory, applications, generalizations},   Cambridge University Press, 2003.              

 

\bibitem{BLS} A.~F.~Beardon, J. Lehner,   and M. Sheingorn,
{\em Closed geodesics on a Riemann surface with application to the Markov spectrum},
Trans. Amer. Math. Soc. 295 (1986), no. 2, 635--647. 
 
\bibitem{Bumby} R.~Bumby, {\em The Markov spectrum} in:  Diophantine approximation and its applications (Proc. Conf., Washington, D.C., 1972), C.~Osgood, ed., pp. 25--58. Academic Press, New York, 1973. 

\bibitem{C} H. Cohn, {\em  Approach to Markoff's minimal forms through
modular functions},
 Ann. of Math. (2)  61 (1955), 1--12.
 
 
\bibitem{CF}  T. W. Cusick and M. E. Flahive, 
{\em The Markoff and Lagrange spectra} (Mathematical Surveys and
Monographs 30, American Mathematical Society, Providence R.I., 1989).  
 
\bibitem{G} 
I. J. Good, {\em  The fractional dimension of continued fractions}, Proc. Camb. Phil. Soc. 37 (1941), 199--228. 

 
\bibitem{H} A. Haas, {\em Diophantine approximation on hyperbolic
Riemann surfaces}, Acta Math. 156 (1986), 33--82.

\bibitem{HS} A. Haas and C. Series, {\em The Hurwitz constant and Diophantine approximation on Hecke groups}, J. London Math. Soc. (2)34 (1986), 219--234.


\bibitem{JP} O. Jenkinson and M. Pollicott, 
{\em Computing the dimension of dynamically defined sets: $E_2$ and bounded continued fractions}, 
Ergodic Theory Dynam. Systems 21 (2001), no. 5, 1429--1445. 

 
\bibitem{Mc} G. McShane, {\em Simple geodesics and a series constant over Teichmuller space},   Invent. Math.  132  (1998),  no. 3, 607--632.

 
\bibitem{vdP} A. van der Poorten, {\em Symmetry and folding of continued fractions}, J. Th\'eor. Nombres Bordeaux 14 (2002), no. 2, 603--611.

 
\bibitem{Sch} A. L. Schmidt, 
{\em Minimum of quadratic forms with respect to Fuchsian groups. I}, J. Reine Angew. Math. 286/287  (1976),  341--368.


\bibitem{scgArt} T. A. Schmidt and M. Sheingorn, 
{\em Parametrizing
simple closed geodesy on
\threesurf}, J. Aust. Math. Soc.  74  (2003),  no. 1, 43--60.

\bibitem{scgMcId} \bysame, 
{\em McShane's identity, using elliptic elements}, Geom. Dedicata (2008) 134:75--90.


\bibitem{ssEG} \bysame, {\em Low height geodesics on \threesurf:    height formulas and examples},  Int. J. Number Theory  3  (2007),  no. 3, 475--501.
 
\bibitem{Sh} M. Sheingorn, {\em Characterization of simple closed geodesics on Fricke surfaces}, Duke Math. J.  52 (1985), 535--545.

 
\bibitem{W} S. Wolpert, {\em  On the K\"ahler form of the moduli space of once punctured tori},  Comment. Math. Helv. 58 (1983), no. 2, 246--256. 


\end{thebibliography}
\end{document}